\documentclass[a4paper,12pt]{amsart}
\usepackage{ifthen}
\usepackage{mathrsfs}
\nonstopmode
\numberwithin{equation}{section}
\newtheorem{thm}{Theorem}

\theoremstyle{definition}

\newenvironment{pf}[1][]{%
 \vskip 3mm
 \noindent
 \ifthenelse{\equal{#1}{}}%
  {{\slshape Proof. }}%
  {{\slshape #1.} }%
 }%
{\qed\bigskip}

\newcounter{alphabet}
\newcounter{tmp}
\newenvironment{Thm}[1][]{\refstepcounter{alphabet}%
\bigskip%
\noindent%
{\bf Theorem \Alph{alphabet}}%
\ifthenelse{\equal{#1}{}}{}{ (#1)}%
{\bf .}
\itshape}{\vskip 8pt}

\newcommand{\C}{{\mathbb C}}

\newcommand{\D}{{\mathbb D}}
\newcommand{\sphere}{{\widehat{\mathbb C}}}

\newcommand{\inv}{^{-1}}

\newcommand{\crit}{{\operatorname{Crit}}}
\newcommand{\poly}{{\mathcal{P}}}
\renewcommand{\setminus}{{-}}

\begin{document}
\bibliographystyle{amsplain}
\title{
Dual mean value problem for complex polynomials
}

\author[V.~Dubinin]{Vladimir Dubinin}
\address{Institute of Applied Mathematics, Far-Eastern Branch, Russian Academy of Sciences,
7 Radio Street, Vladivostok, 690041, Russia}
\email{dubinin@iam.dvo.ru}
\author[T. Sugawa]{Toshiyuki Sugawa}
\address{Graduate School of Information Sciences,
Tohoku University, Aoba-ku, Sendai 980-8579, Japan}
\email{sugawa@math.is.tohoku.ac.jp}
\keywords{Smale's mean value conjecture, critical point}
\subjclass[2000]{Primary 30C10; Secondary 30C55}
\date{June 22, 2009,\quad {\tt File: dubinin-sugawa0906.tex}}
\begin{abstract}
We consider an extremal problem for polynomials, which is dual to the well-known
Smale mean value problem.
We give a rough estimate depending only on the degree.
\end{abstract}
\thanks{
The present research was partially supported by 
the Russian Foundation for Basic Research (grant no.~08-01-00028),
the Far-Eastern Branch of RAS (grant no.~09-I-P4-02),
%the Far-Eastern Branch of RAS (grant no.~06-III-A-01-013),
and the JSPS Grant-in-Aid for Scientific Research (B), 17340039.}
\maketitle

\section{Introduction}
Let $P$ be a complex polynomial of degree $d\ge2,$ that is,
$P$ is a polynomial map on the complex plane $\C$ of the form
$$
P(z)=c_dz^d+c_{d-1}z^{d-1}+\dots+c_1z+c_0
$$ 
for complex coefficients $c_0,\dots, c_d$ with $c_d\ne0.$
We denote by $\poly_d$ the set of complex polynomials of degree $d.$
We say that $\zeta\in\C$ is a {\it critical point} of $P$
if $P'(\zeta)=0.$
The image $P(\zeta)$ of a critical point $\zeta$ under $P$ is called
a {\it critical value}.
We denote by $\crit(P)$ the set of critical points of $P.$
In 1981, Smale \cite{Smale81} proved the following inequality
in connection with root-finding algorithms.

\begin{Thm}[Smale]
Let $P$ be a polynomial of degree $d\ge2$ over $\C$ and
suppose that $z\in\C$ is not a critical point of $P.$
Then there exists a critical point $\zeta$ of $P$ such that
\begin{equation}\label{eq:S}
\left|
\frac{P(\zeta)-P(z)}{\zeta-z}
\right|
\le 4|P'(z)|.
\end{equation}
\end{Thm}

In the same paper, Smale asked whether the factor $4$ can be replaced
by $1$ or even by $1-1/d.$
See also \cite{SS86} and \cite[\S~7.2]{RS:poly} 
for background and further references.

To simplify expressions and to emphasize invariance, we introduce some notation.
For $\zeta\in\crit(P)$ and $z\in\C\setminus\crit(P),$ we define $Q(P,z,\zeta)$ by
$$
Q(P,z,\zeta)=\frac{P(z)-P(\zeta)}{(z-\zeta)P'(z)}.
$$
%when $z\notin\crit(P)$ and by $Q(P,z,\zeta)=\lim_{z'\to z}Q(P,z',\zeta)$
%when $z\in\crit(P).$
%It is easy to see that $Q(P,z,\zeta)=\infty$ for 
%$z\in\crit(P)\setminus\{\zeta\}$
%and $Q(P,\zeta,\zeta)=1/k,$ where $k\ge2$ is the order of $P-P(\zeta)$
%at $\zeta:$
%$P'(\zeta)=\dots=P^{(k-1)}(\zeta)=0\ne P^{(k)}(\zeta).$
It is easy to verify the invariance relation
$Q(\tilde P,\tilde z,\tilde\zeta)=Q(P,z,\zeta)$
for $z=a\tilde z+b, \zeta=a\tilde\zeta+b$ and
$\tilde P(\tilde z)=AP(a\tilde z+b)+B$ for constants $a,b,A,B$ with
$aA\ne0.$

We further set
$$
S(P,z)=\min
\{|Q(P,z,\zeta)|
%\left\{\left|\frac{P(\zeta)-P(z)}{(\zeta-z)P'(z)}\right|
: \zeta\in\crit(P)\}
%\right\}
$$
for $z\in\C$ and
$$
K(d)=\sup\{S(P,z): P\in\poly_d, z\in\C\}.
$$
Smale's theorem says that $K(d)\le 4,$ while the example
$P_0(z)=z^d-dz$ shows that $K(d)\ge 1-1/d.$
Smale's problem asks to find the value of $K(d)$
and Smale's conjecture can be stated as $K(d)=1-1/d.$
This conjecture has been confirmed for degrees $d=2,3,4$ 
(Tischler \cite{Tis89}), 
for $d=5$ (Crane \cite{CraneCP}). 
Sendov and Marinov \cite{SM07} claim that the conjecture is
true for $d\le 10$ by massive numerical computations, which
one cannot check easily.
Meanwhile, some improvements were made for Smale's theorem,
see \cite{BMN02}, \cite{CFL07}, \cite{FS06}, \cite{Crane07}.
Note that all the known estimates $K(d)\le \tilde K(d)$ satisfy
$\liminf_{d\to\infty} \tilde K(d)\ge4.$

We may pose a dual problem to Smale's mean value problem:
Consider the quantity
$$
T(P,z)=\max\{|Q(P,z,\zeta)|: \zeta\in\crit(P)\}
$$
for $z\in\C$ and find the value
$$
L(d)=\inf\{T(P,z): P\in\poly_d,  z\in\C\}.
$$
For the polynomial $P^*(z)=(z+1)^d-1,$ we have $T(P^*,0)=1/d,$ and
hence $L(d)\le 1/d.$
%It is likely that $L(d)=1/d$ for $d\ge2.$

We should now mention Tischler's strong form of Smale's conjecture:
\begin{equation}\label{eq:T}
\min\left\{\left|
Q(P,z,\zeta)-\frac12
\right|: \zeta\in\crit(P)\right\}\le \frac12-\frac1d,
\quad P\in\poly_d.
\end{equation}
By the triangle inequality $||Q(P,z,\zeta)|-1/2|\le |Q(P,z,\zeta)-1/2|,$
\eqref{eq:T} would imply the inequality:
$$
\max\left\{S(P,z)-\frac12, \frac12-T(P,z)\right\}
\le \frac12-\frac1d,
$$
in other words, $S(P,z)\le 1-1/d$ and $1/d\le T(P,z).$
Since the inequalities $1-1/d\le K(d)$ and $L(d)\le 1/d$ trivially hold,
the relations $K(d)=1-1/d$ and $L(d)=1/d$ would follow from \eqref{eq:T}.

It is known that \eqref{eq:T} is valid for $d\le4,$ see
Tischler \cite{Tis89}, where the case $d=4$ is attributed to J.-C. Sikorav.
In particular, we have $L(d)=1/d$ for $d=2,3,4.$
Unfortunately, \eqref{eq:T} does not hold in general.
Indeed, Tyson \cite{Tyson05} revealed that \eqref{eq:T} does not hold
for every $d\ge 5.$
%Note that Sendov and Marinov \cite{SM07} numerically found for $d\le 10,$
%the best constants $\tau_d\in[0,1/2]$ such that the inequality
%$$
%\min\left\{\left|
%Q(P,z,\zeta)-\frac12+\tau_d
%\right|: \zeta\in\crit(P)\right\}\le \tau_d+\frac12-\frac1d
%$$
%holds for $P\in \poly_d, z\in\C, \zeta\in\crit(P)$
%to prove $K(d)=1-1/d$ for $d\le10.$
%Note that their result simultaneously implies $L(d)\ge 1/d-2\tau_d.$
%For instance, $L(10)\ge 1/10-0.0003178\dots$ by their computations.

In the present note, we give the rough estimate $L(d)\ge 1/(d4^d)$
as the first step towards the conjecture $L(d)=1/d.$

\begin{thm}\label{thm:main1}
Let $P$ be a polynomial of degree $d\ge2$ over $\C$ and
suppose that $z\in\C$ is not a critical point of $P.$
Then there exists a critical point $\zeta$ of $P$ such that
\begin{equation*}\label{eq:main}
\frac{|P'(z)|}{d4^d}\le
\left|
\frac{P(\zeta)-P(z)}{\zeta-z}
\right|.
\end{equation*}
\end{thm}

We proved this theorem when the first author visited Hiroshima
University, where the second author was working, in July 2007.
We learnt from T.~W.~Ng that he considered the same problem independently
and that he showed the inequality $L(d)>0$ for each $d$ in the same year
by using the notion of amoebae (see \cite{Ng07}).
Since there is no explicit lower bound for $L(d)$ in the literature so far,
our result seems to be meaningful even though the bound is too far from the
conjectured one.

Note also that the conjecture is true for the special case 
when $P$ is conservative, 
that is, $P(\zeta)=\zeta$ for every $\zeta\in\crit(P).$
See the last corollary in \cite[p.~455]{Tis89}, 
where Tischler attributes its parent theorem to Yoccoz.

\section{Proof of the theorem}

By the transformation $\tilde P(\tilde z)=[P(\tilde z+z_0)-P(z_0)]/c_d,$
one can easily see that Theorem \ref{thm:main1}
is equivalent to the following assertion.

\begin{thm}\label{thm:main}
Let $\zeta_1,\dots,\zeta_{d-1}$ be the critical points of a monic polynomial
$P$ of degree $d\ge2$ with $P'(0)\ne0.$
Then
$$
\max_{1\le j\le d-1}\left|\frac{P(\zeta_j)}{\zeta_j P'(0)}\right|
\ge \frac1{d4^d}.
$$
\end{thm}

\begin{pf}
%Let $P(z)=z^d+c_{d-1}z^{d-1}+\dots+c_1z.$
Since $P'(z)=d(z-\zeta_1)\cdots(z-\zeta_{d-1}),$ the relation
\begin{equation}\label{eq:c1}
P'(0)=d(-1)^{d-1}\zeta_1\cdots\zeta_{d-1}
\end{equation}
holds.

For $t>0,$ we set $\Delta(t)=\{w\in\sphere: |w|>t\},$
where $\sphere$ stands for the Riemann sphere $\C\cup\{\infty\}.$
Let
$$
R=\max_{1\le j\le d-1}|P(\zeta_j)|^{1/d}
$$
and set $A=\{z\in\C: |P(z)|\le R^d\}.$
Note that $A$ is a continuum containing $0$ and all the critical points $\zeta_1,
\dots, \zeta_{d-1}.$
%Let $A$ be the union of the segments $[0,\zeta_j],~ j=1,\dots, d-1$ and
%$B=\{w\in\C: w^d\in A\}.$
Since $P:\C-A\to\Delta(R^d)-\{\infty\}$ is an unbranched covering map of degree $d,$
one can take a single-valued analytic branch $F(w)$ of
$P\inv(w^d)$ on $\Delta(R)$ with $F(w)=w+O(1)$ as $w\to\infty.$
Note that the function $F$ satisfies the relation
\begin{equation}\label{eq:fe}
P(F(w))=w^d,\quad |w|>R
\end{equation}
and that $F(\Delta(R))=\sphere-A.$
We now show that $F$ is univalent in $\Delta(R).$
Suppose, to the contrary, that $F(w_0)=F(w_1)$
for distinct points $w_0$ and $w_1$ in $\Delta(R).$
Obviously, $z_0=F(w_0)$ is not the point at infinity.
Since $P$ has no critical point in $\C-A,$
$P$ is univalent in a small neighborhood $V$ of $z_0.$
On the other hand, by \eqref{eq:fe}, we have $w_0^d=w_1^d,$
and thus, $w_1=\tau w_0$ for a complex number $\tau\ne1$ with
$\tau^d=1.$
We consider the function $F_1(w)=F(\tau w)$ in $|w|>R.$
Then $F(w_0)=F_1(w_0).$
Since $F_1(w)=\tau w+O(1)$ as $w\to\infty,$
$F_1$ and $F$ are not identically equal on $\Delta(R).$
Therefore, $F(w)\ne F_1(w)$ for $w\ne w_0$ close enough to $w_0.$
By \eqref{eq:fe}, we have
$P(F(w))=P(F_1(w))$ for $w\in\Delta(R).$
This is impossible because $P$ is univalent in $V.$
Thus we have shown that $F$ is univalent in $\Delta(R).$

Let $f$ be a function $f$ on the unit disk
$\D=\{z\in\C: |z|<1\}$ by $f(z)=R/F(R/z).$
Since $F$ has no zero in $\Delta(R),$ we see that
$f$ is analytic in $\D$ and that $f(0)=f'(0)-1=0.$
From what we saw above, we also obtain that $f$ is univalent
in $\D.$
The Koebe one-quarter theorem implies that
the image $f(\D)$ contains the disk $|w|<1/4.$
Since $\zeta_j\in A=\sphere\setminus F(\Delta(R)),$ we have
\begin{equation}\label{eq:R}
\frac R{|\zeta_j|}\ge \frac14,\quad j=1,2,\dots, d-1.
\end{equation}
We may assume that $R=|P(\zeta_1)|^{1/d}.$
Then, by \eqref{eq:c1} and \eqref{eq:R}, we have
$$
\max_{1\le j\le d-1}\left|\frac{P(\zeta_j)}{\zeta_j P'(0)}\right|
\ge \left|\frac{P(\zeta_1)}{\zeta_1 P'(0)}\right|
=\frac{R^d}{|\zeta_1|\cdot d|\zeta_1\cdots\zeta_{d-1}|}
\ge \frac{1}{d4^d}.
$$
The proof is now complete.
\end{pf}

\def\cprime{$'$} \def\cprime{$'$} \def\cprime{$'$}
\providecommand{\bysame}{\leavevmode\hbox to3em{\hrulefill}\thinspace}
\providecommand{\MR}{\relax\ifhmode\unskip\space\fi MR }
% \MRhref is called by the amsart/book/proc definition of \MR.
\providecommand{\MRhref}[2]{%
  \href{http://www.ams.org/mathscinet-getitem?mr=#1}{#2}
}
\providecommand{\href}[2]{#2}

%\bibliography{papers}
\end{document}